\documentclass[11pt]{amsart}
\usepackage{amssymb,amsmath}
\usepackage{diagbox}
\usepackage{pgfplots}
\pgfplotsset{compat = newest}

\usepackage[pagebackref]{hyperref} 

\newcommand{\QQ}{{\mathbb{Q}}}
\newcommand{\RR}{{\mathbb{R}}}
\newcommand{\ZZ}{{\mathbb{Z}}}

\newcommand{\fA}{{\mathfrak{A}}}
\newcommand{\fS}{{\mathfrak{S}}}

\newcommand{\Aut}{{\operatorname{Aut}}}
\newcommand{\Gal}{{\operatorname{Gal}}}
\newcommand{\PGL}{{\operatorname{PGL}}}
\newcommand{\PSL}{{\operatorname{PSL}}}
\newcommand{\pr}{{\mathrm {pr}}}
\newcommand{\E}{{\operatorname{E-}}}
\newcommand{\Ex}{{\operatorname{E}}}
\newcommand\gp{{\tt pari-gp}}

\let\eps=\epsilon
\let\cd=\cdot
\let\ti=\times

\newtheorem{thm}{Theorem}[section]
\newtheorem{prop}[thm]{Proposition}

\begin{document}
\title{Computational data on $\fS_n$-extensions of $\QQ$}

\author{Gunter Malle}
\address{FB Mathematik, RPTU Kaiserslautern,
  Postfach 3049, D--67653 Kaisers\-lautern, Germany.}
\makeatletter
\email{malle@mathematik.uni-kl.de}
\makeatother

\date{\today}

\thanks{The author gratefully acknowledges support by the Deutsche
 Forschungsgemeinschaft --- Project-ID 286237555 -- TRR 195.}

\begin{abstract}
We discuss computational results on field extensions $K/\QQ$ of
degree~$n\le11$ with Galois group of the Galois closure isomorphic to the
full symmetric group $\fS_n$. More precisely, we present statistics on the
number of such extensions as a function of the field discriminant and compare
them to the known predictions by Bhargava and the author. We also investigate
the numbers of fields with equal discriminant and tabulate class numbers and
class groups to compare them against Cohen--Lenstra--Martinet type of
heuristics and their proposed improvements.
\end{abstract}

\maketitle


In this note we report on extensive computer calculations pertaining to field
extensions $K/\QQ$ (in a fixed algebraic closure) of degree $n\le11$ with
Galois group of the Galois closure isomorphic to the full symmetric group
$\fS_n$. We present statistics on the growth behaviour of the counting function
for such fields, ordered by (absolute) discriminant and compare them to
conjectural predictions by Bhargava \cite{Bh07} and the author \cite{Ma02}.
Assuming these predictions hold, our data also give some indication on the
magnitude of the error term. We compare the
distribution of class groups for these fields with the predictions made by
Cohen--Lenstra and Cohen--Martinet \cite{CM87}, and more recently by the author
\cite{Ma08,Ma10} and by Sawin--Wood \cite{SW23} in the presence of
roots of unity. Our most extensive results are for fields with at most one real
embedding, but we also record data on totally real extensions.  As a byproduct
we determine the minimal discriminants for several pairs of Galois groups and
signatures for which these were previously unknown.

\section{Enumeration of $\fS_n$-extensions with at most one real embedding}

In this section we display data on $\fS_n$-extensions of $\QQ$ of bounded
discriminant, for $4\le n\le 11$. We compare these to the predictions
coming from Bhargava's conjecture \cite[Conj.~1.2]{Bh07}.

\subsection{Growth rate and Bhargava's constant}   \label{subsec:growth}
For an integer $n\ge2$ and a real number $X\ge0$ we let
$$N_t(\fS_n,X):=\#\{K/\QQ : |D(K/\QQ)|\le X,\ \Gal(K/\QQ)\cong\fS_n\}$$
be the number of degree~$n$ extensions $K$ of $\QQ$ (inside a fixed algebraic
closure) of absolute discriminant bounded by $X$ and with Galois group of the
Galois closure the full symmetric group $\fS_n$ (here denoted $\Gal(K/\QQ)$ by
abuse of notation). By a result of Schmidt, the total number of
degree~$n$-extensions of $\QQ$ of absolute discriminant at most $X$ is
$O(X^{(n+2)/4})$; the currently best known upper bound for arbitrary $n\ge6$ is
$O(X^{1.564(\log n)^2})$, shown by Lemke Oliver and Thorne \cite{LOT20} (see
also Couveignes \cite{Co20}, and the recent far reaching improvements for large
$n$ by Lemke Oliver \cite{LO23}). We conjectured in \cite[p.~2]{Ma02} that
asymptotically, $N_t(\fS_n,X)$ grows linearly with $X$, that is,
$$\lim_{X\to\infty} N_t(\fS_n,X)/X=c_n$$
for a positive constant $c_n$. Bhargava \cite[Conj.~1.2]{Bh07} proposes an
explicit Euler product for the value of~$c_n$. Note that this conjecture
is known to hold when $n=2,3$ by classical results, so we don't include these
cases in our table. Also, Bhargava \cite{Bh05,Bh10} proved his conjecture for
$n=4,5$, and our data fit quite nicely to this. See also \cite{ST23} for a
recent further heuristic explanation of Bhargava conjecture.

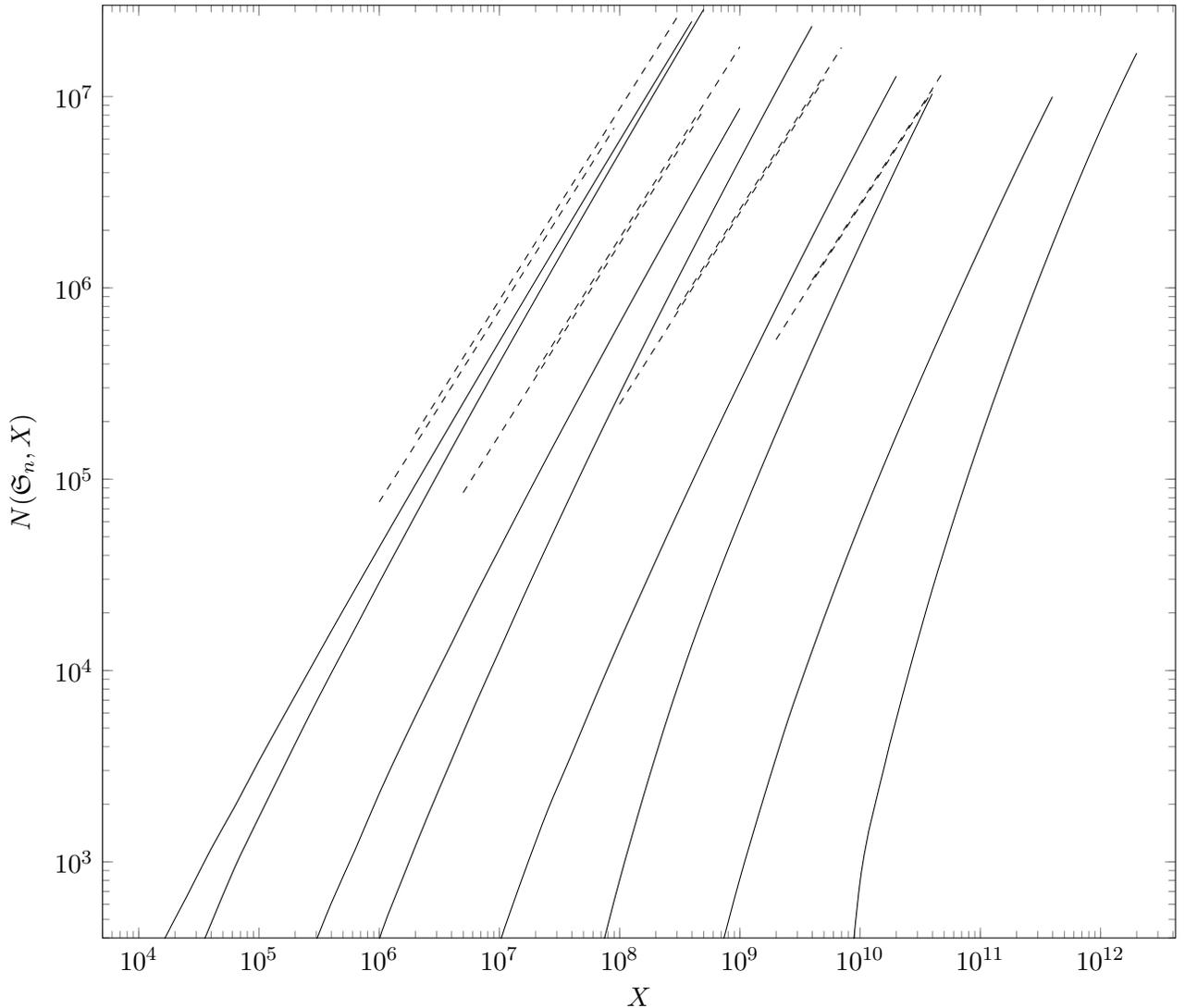
\begin{table}[htbp]
\caption{Number of $\fS_n$-extension with $r_1\le1$, for $4\le n\le 11$}   \label{tab:plot}
\begin{tikzpicture}
\begin{loglogaxis}[
    legend style = { at = {(0.25,0.95)}},
    xmin = 5000, xmax = 4.2*10^12,
    ymin = 400, ymax = 3*10^7,
    xtick distance = 10,
    ytick distance = 10,
    grid = none,
    minor tick num = 1,
    width = \textwidth,
    height = 15cm, 
    xlabel = {$X$},
    ylabel = {$N(\fS_n,X)$},]
 
\addplot[domain=10^6:9*10^7, dashed]{0.076043142*x}; 
\addplot[ thin, smooth ]
coordinates {
(6309.573,119)
(10000.00,206)
(15848.931,380)
(25118.864,661)
(39810.717,1165)
(63095.734,1954)
(100000.00,3374)
(158489.319,5691)
(251188.643,9550)
(398107.170,16010)
(630957.344,26637)
(1000000.00,44122)
(1584893.192,72900)
(2511886.432,119885)
(3981071.706,196460)
(6309573.445,321736)
(10000000.00,525099)
(15848931.92461,854975)
(25118864.31510,1388766)
(39810717.05535,2251820)
(63095734.44802,3646175)
(100000000.,5892088)
(158489319.,9510689)
(251188643.,15330201)
(398107170.,24678890) };  
\addplot[domain=2*10^6:3*10^8, dashed]{0.08635053*x}; 
\addplot[ thin, smooth ]
coordinates {
(15848.931,121)
(25118.864,237)
(39810.717,474)
(63095.734,940)
(100000.00,1714)
(158489.3192,3079)
(251188.6432,5515)
(398107.1706,9654)
(630957.3445,16620)
(1000000.00,28831)
(1584893.19,49403)
(2511886.43,84346)
(3981071.70,143466)
(6309573.44,242492)
(10000000.0,406505)
(15848931.9,678359)
(25118864.3,1128225)
(39810717.0,1870684)
(63095734.4,3090861)
(100000000.0,5089844)
(158489319.2,8363098)
(251188643.1,13703390)
(398107170.5,22408460)
(500000000.0,28552958) };  
\addplot[domain=5*10^6:5*10^8, dashed]{0.01702530105*x}; 
\addplot[ thin, smooth ]
coordinates {
(158489.3192,153)
(251188.6432,296)
(398107.1706,605)
(630957.3445,1158)
(1000000.000,2276)
(1584893.192,4247)
(2511886.432,7710)
(3981071.706,13764)
(6309573.445,24698)
(10000000.00,43254)
(15848931.92,75622)
(25118864.32,131151)
(39810717.06,225033)
(63095734.45,384682)
(100000000.0,654757)
(158489319.2,1107168)
(251188643.1,1864703)
(398107170.5,3127652)
(630957344.4,5220176)
(1000000000,8681205) };  
\addplot[domain=2*10^7:10^9, dashed]{0.0182218533*x }; 
\addplot[ thin, smooth ]
coordinates {
(630957.344,129)
(1000000.00,392)
(1584893.19,837)
(2511886.43,1715)
(3981071.70,3386)
(6309573.44,6673)
(10000000.0,12698)
(15848931.9,24283)
(25118864.3,45439)
(39810717.0,84209)
(63095734.4,154538)
(100000000.0,279385)
(158489319.2,498803)
(251188643.1,882449)
(398107170.5,1549564)
(630957344.4,2704109)
(1000000000.,4690185)
(1584893192.,8076149)
(2511886431.,13816562)
(3981071705.,23330805) };  
\addplot[domain=10^8:5*10^9, dashed]{0.002468803218*x }; 
\addplot[ thin, smooth ]
coordinates {
(6309573.4448,138)
(10000000.000,373)
(15848931.925,855)
(25118864.315,1850)
(39810717.055,3648)
(63095734.448,7293)
(100000000.00,14261)
(158489319.25,27157)
(251188643.15,51492)
(398107170.55,95916)
(630957344.48,177048)
(1000000000.0,322797)
(1584893192.4,582440)
(2511886431.5,1039331)
(3981071705.5,1838540)
(6309573444.8,3227601)
(10000000000.0,5627311)
(15848931924.6,9730933)
(20000000000,12778743) };  
\addplot[domain=3*10^8:7*10^9, dashed]{0.0025736805608283*x }; 
\addplot[ thin, smooth ]
coordinates {
(39810717.055,64)
(63095734.448,257)
(100000000.0,803)
(158489319.2,2225)
(251188643.1,5707)
(398107170.5,13487)
(630957344.4,29472)
(1000000000.0,61353)
(1584893192.5,123268)
(2511886431.5,242217)
(3981071705.5,468298)
(6309573444.8,892039)
(10000000000.0,1676193)
(15848931924.6,3114931)
(25118864315.1,5717440)
(39810717055.3,10322101) };  
\addplot[domain=2*10^9:4*10^10, dashed]{0.0002683962162*x }; 
\addplot[ thin, smooth ]
coordinates {
(251188643.15,11)
(398107170.55,83)
(630957344.48,274)
(1000000000.0,817)
(1584893192.4,2181)
(2511886431.5,5509)
(3981071705.5,12625)
(6309573444.8,27777)
(10000000000.0,58671)
(15848931924.6,119264)
(25118864315.1,235327)
(39810717055.3,456134)
(63095734448.0,869337)
(100000000000,1634855)
(158489319246,3044442)
(251188643151,5573423)
(398107170553,9972448) };  
\addplot[domain=4*10^9:5*10^10, dashed]{0.0002747762886*x }; 
\addplot[ thin, smooth ]
coordinates {
(6309573444.802,38)
(10000000000.00,789)
(15848931924.61,3113)
(25118864315.10,9518)
(39810717055.35,26489)
(63095734448.02,67625)
(100000000000.0,162093)
(158489319246.1,368111)
(251188643151.0,801689)
(398107170553.5,1682873)
(630957344480.2,3417569)
(1000000000000,6698503)
(1584893192461,12537148)
(2000000000000,16863137) };  
\end{loglogaxis}
 
\end{tikzpicture}
\end{table}

In Table~\ref{tab:plot} we display our computational results on the number
$N(\fS_n,X)$ of $\fS_n$-extensions with bounded absolute discriminant $X$ and
with at most one real embedding. More precisely, we give eight curves for
$N(\fS_n,X)$, where $n=4,\ldots,11$ from left to right, in doubly logarithmic
scale. Also, as dashed lines we indicate the predicted (linear) asymptotic
growth $\sim B_n\,X$ according to Bhargava. The underlying numerical data will
be presented in the subsequent sections.

Approximations of the relevant Bhargava constants $B_n$ for $\fS_n$-extensions
with at most one real embedding are given in Table~\ref{tab:Bhar}. Note that
the constants for even degree $n=2m$ are slightly smaller than those for the
next larger odd degree $n=2m+1$, in our range $m\le 5$. Thus, asymptotically we
expect more $\fS_{2m+1}$-fields than $\fS_{2m}$-fields up to a given
discriminant bound!

\begin{table}[htb]
\caption{Bhargava constants $B_n$ for $r_1\le1$}   \label{tab:Bhar}
$\begin{array}{cc|cc}
 n& B_n& n& B_n\\
\hline
  4& 0.07604314\ldots&  5& 0.08635053\ldots\\
  6& 0.01702530\ldots&  7& 0.01822185\ldots\\
  8& 0.00246880\ldots&  9& 0.00257368\ldots\\
 10& 0.00026840\ldots& 11& 0.00027478\ldots\\
\hline
\end{array}$
\end{table}

The basic algorithm for obtaining complete tables of number fields with given
non-solvable Galois group is still Hunter's method, explained for example in
\cite[\S9.3.1]{Coh}. It uses the fact that a primitive number field with
discriminant $D$ contains a primitive integral element whose minimal polynomial
has coefficients that can be explicitly bounded in terms of $D$. Hence all such
fields can be found by enumerating all monic irreducible polynomials with
integer coefficients satisfying these coefficient bounds and then discarding
those with larger field discriminant, wrong signature or Galois group, and
eliminating isomorphic fields. Our results were obtained by implementing this
method in the \gp{} system \cite{pari}. We used the command {\tt polredabs} to
find generating polynomials with small coefficients. In the end we tested all
pairs of fields with same discriminant on isomorphism: no two distinct reduced
polynomials turned out to define isomorphic fields.

Unfortunately this method is not very efficient since exponential in the degree
and moreover tends to produce many isomorphic fields as well as fields of
larger discriminant. This limits
the applicability severely. For degree $n\ge7$, our computations did not
complete but heuristic reasoning indicates that the obtained numbers of fields
should be within $0.01\%$ of the correct ones at least for $n\le9$.

\section{Totally complex $\fS_4$-extensions}   \label{sec:S4}
Using Hunter's method we have enumerated the primitive totally complex number
fields of degree~4 (up to isomorphism) with discriminant at most $10^9$.

\subsection{The fields}
The only possible Galois groups (of the Galois closure) are the symmetric and
the alternating group on four letters. We note that the degree~4 fields with
alternating Galois group $\fA_4$ and discriminant up to $10^{13}$ had already
been counted in \cite[\S8]{CDO}.

\begin{prop}   \label{prop:deg4}
 There are $63\,748\,067$ totally complex number fields of degree~$4$ (up to
 isomorphism) with discriminant at most $10^9$ and Galois group $\fS_4$.
\end{prop}

Our data for $\fS_4$-extensions are summarised in Table~\ref{tab:S4}. In
\cite[Tab.~9.3]{CDO} the authors also gave the number of totally complex
$\fS_4$-fields of discriminant at most~$10^7$ (which agrees with ours).
Note that the prediction for totally complex $\fS_4$-fields up to $10^8$ made
in \cite[9.2]{CDO} is within 0.5\% of the correct value.

\begin{table}[ht]
\caption{Totally complex $\fS_4$-fields}  \label{tab:S4}
$\begin{array}{|c|r|r|r|r|r|r|}
\hline
    X& N(\fS_4,X)& E(\fS_4,X)& \qquad\alpha_1& \qquad\alpha_2& P(\fS_4,X)& \Delta(\fS_4,X)\cr
\hline\hline
        10^3&        8&       68&  .611& .918&        9& .00\cr
        10^4&      206&      554&  .686& .895&      216& .24\cr
        10^5&     3374&     4230&  .725& .887&     3377& .06\cr
        10^6&    44122&    31921&  .751& .873&    44163& .27\cr
        10^7&   525099&   235332&  .767& .865&   525012& .28\cr
        10^8&  5892088&  1712226&  .779& .860&  5892769& .35\cr
        10^9& 63748067& 12295075&  .788& .855& 63748073& .09\cr
\hline
\end{array}$   
\end{table}

By the theorem of Bhargava \cite{Bh05} the number of totally complex
$\fS_4$-fields asymptotically grows linearly with the discriminant, with
proportionality
constant $B_4$ (see Table~\ref{tab:Bhar}). It seems interesting to obtain
information on the error term in this asymptotic behaviour. For this
assume that
$$E(\fS_4,X):=B_4X-N(\fS_4,X)\sim c' X^\alpha\qquad\text{(for $X\to\infty$)}$$
for some $\alpha<1$ and $c'\in\RR^\times$. Then
$$\alpha_1(\fS_4,X):=\log E(\fS_4,X)/\log X$$
and
$$\alpha_2(\fS_4,X):=\log \big(E(\fS_4,X)/E(\fS_4,X/2)\big)/\log(2)$$
both converge to $\alpha$ when $X\rightarrow\infty$. In
the fourth and fifth column of Table~\ref{tab:S4} we give
$\alpha_1(\fS_4,X)$ and $\alpha_2(\fS_4,X)$. In the range of our data, the
first increases monotonically, while the second decreases, suggesting
that the exponent $\alpha$ of the error term should satisfy
$0.788 \le \alpha\le 0.855$.   
This is consistent with the results in \cite{Ma06} where we enumerated the
totally real $\fS_4$-fields of degree~4 up to discriminant $10^9$ and found
that the correspondingly defined $\alpha$ should lie between $0.76$ and $0.87$.
\par
The authors of \cite{CDO} speculate, following a communication by Yukie,
that $\alpha=5/6$ (which would fit with the above observations), and more
precisely an asymptotic behaviour of the form
$$N(\fS_4,X) = B_4 X + c' X^{\frac{5}{6}} + c'' X^{\frac{3}{4}}\log X
  + c''' X^{\frac{3}{4}} + O(X^\beta)\eqno{(1)}$$
with some exponent $\beta<3/4$. Using their count up to discriminant $10^7$ they
present least squares approximations to the constants $c',c'',c'''$ (which
then predict 5\,902\,307 fields up to $10^8$, less than 0.2\% off the actual
value).
Assuming this form of the asymptotic expansion, using a least square method
with evaluation points at $X=10^{i/10}$, $i\ge20$, on our more extensive
results we get the approximations
$$c'=-.50856,\qquad c''=.0106,\qquad c'''=.4533,$$    
to the constants appearing in the conjectured expansion~$(1)$ for $N(\fS_4,X)$.
The value of
$$P(\fS_4,X):=B_4 X + c' X^{\frac{5}{6}} + c'' X^{\frac{3}{4}}\log X
  + c''' X^{\frac{3}{4}},$$
with $c', c'',c'''$ as before, is given in column~6 of Table~\ref{tab:S4}.
In the last column we display the quantity
$$\Delta(\fS_4,X):=\log|N(\fS_4,X)-P(\fS_4,X)|/\log(X),$$
which is consistent with $\beta<3/4$ in~$(1)$.
The formula predicts~672\,934\,742 totally complex $\fS_4$-fields of
discriminant at most~$10^{10}$.

\subsection{Discriminant multiplicities} We have also counted the occurrences of
$\fS_4$-fields with the same discriminant:

\begin{prop}   \label{prop:mult4}
 The multiplicities $k\le24$ of discriminants of the first $N$ totally complex
 $\fS_4$-fields of degree~$4$, where $N\in\{5\cdot 10^6,10^7,5\cdot 10^7\}$,
 are as given in Table~$\ref{tab:discs4}$.
\end{prop}

\begin{table}[ht]
\caption{$k$-tuples of $\fS_4$-discriminants}  \label{tab:discs4}
$\small\begin{array}{|r|r|r|r|r|r|r|r|r|r|r|r|r|r|r|}
\hline
   \quad k&  1& 2& 3& 4& 5& 6& 7& 8& 9& 10& 11& 12\cr
    N\quad\ & 13& 14& 15& 16& 17& 18& 19& 20& 21& 22& 23& 24\cr
\hline
 5\cdot 10^6& 3057365& 673635& 83788& 51683& 14702& 5744& 1685& 1031& 318& 193& 86& 162\\
     & 39& 30& 15& 15& 5& 2& 2& .& 1& .& .& .\\  
 10^7& 5992096& 1374259& 170145& 111445& 31722& 12612& 3756& 2435& 767& 497& 216& 376\\
     & 114& 75& 44& 33& 15& 9& 6& 2& 2& .& .& 1\\  
 5\cdot 10^7& 28791743& 7098100& 867993& 645901& 179889& 74674& 23010& 17494& 5317& 4034& 1442& 2798\\
  & 838& 681& 396& 310& 134& 96& 64& 27& 18& 16& 1& 11\\
\hline
\end{array}$   
\end{table}

There are further 7 discriminants with multiplicities $25,26,28,31,36,37,42$
among the first $5\cdot 10^7$ fields,
the most frequent one being $51\,4764\,864=2^6.3^2.191.4679$.
Kl\"uners \cite{Kl06} has shown that the number of $\fS_4$-fields with
discriminant~$D$ is at most equal to $O_\eps(D^{1/2+\eps})$ for all $\eps>0$.
Conjecturally, it should be much smaller, rather of the order $O_\eps(D^\eps)$,
an expectation consistent with our data in Table~\ref{tab:discs4}.

\subsection{Class groups statistics}
We have calculated the class groups $H_K$ of all $\fS_4$-fields $K$ in our
range using the {\tt bnfinit}-command in \gp{} (so the correctness of the
results relies on a heuristic strengthening of the generalised Riemann
hypothesis, as described in the {\tt pari}-manual).

\begin{prop}   \label{prop:clno4}
 The distribution of odd parts $h'<100$ of the class number among the first
 $6\cd10^7$ totally complex $\fS_4$-fields is displayed in
 Table~$\ref{tab:clno4}$; there are another $120\,040$ such fields with
 $h'>100$.
\end{prop}

\begin{table}[htb]
\caption{Odd parts of class numbers of the first $6\cd10^7$ totally complex
  $\fS_4$-fields}  \label{tab:clno4}
$\begin{array}{|c|r||c|r||c|r||c|r||c|r|}
\hline
    h'& n(h')&  h'& n(h')&  h'& n(h')&  h'& n(h')&  h'& n(h')\\
\hline
  1& 46436885& 21& 154675& 41& 25066& 61& 10345& 81& 8074\\
  3&  6953434& 23&  86703& 43& 22456& 63& 16043& 83& 5053\\
  5&  2221638& 25&  86730& 45& 35089& 65& 11334& 85& 6174\\
  7&  1074329& 27&  89799& 47& 18559& 67&  8195& 87& 6261\\
  9&   828276& 29&  52404& 49& 19092& 69& 10835& 89& 4268\\
 11&   412677& 31&  45455& 51& 21840& 71&  7269& 91& 4990\\
 13&   290462& 33&  57180& 53& 14231& 73&  6770& 93& 5335\\
 15&   325691& 35&  46731& 55& 16870& 75& 10840& 95& 4723\\
 17&   165224& 37&  31157& 57& 16890& 77&  7382& 97& 3543\\
 19&   130751& 39&  39912& 59& 11227& 79&  5758& 99& 5335\\
\hline\end{array}$
\end{table}

We next compare the observed relative frequency of odd parts of class groups
with the heuristic of Cohen--Martinet \cite[(8.1)]{CM87}. They predict that the
relative proportions of odd parts of class groups for totally complex quartic
$\fS_4$-fields should agree with their relative proportions for real quadratic
fields. In Table~\ref{tab:clno4rel} we list the relative proportions for the
smallest $\fS_4$-fields, grouped in packages of 20 million, while in the last
line is printed the Cohen--Martinet prediction. Here, $3^2$ denotes a class
group $(\ZZ/3\ZZ)^2$.

\begin{table}[htb]
\caption{Relative proportions of odd parts $H'$ of $\fS_4$-class groups}
  \label{tab:clno4rel}
$\begin{array}{|r|r|r|r|r|r|r|r|r|r|r|r|r|r|r|}
\hline
 \hbox{\diagbox{$N$}{$H'$}}& 1& 3& 5& 7& 9& 3^2& 11& 13& 15\cr
\hline
 0$--$2\!\cd\! 10^7&  .7799& .114& .366\E1& .178\E1& .122\E1& .105\E2& .680\E2& .477\E2& .515\E2\\
 2$--$4\cd 10^7&  .7720& .117& .371\E1& .180\E1& .128\E1& .121\E2& .689\E2& .488\E2& .551\E2\\
 4$--$6\cd 10^7&  .7699& .117& .373\E1& .179\E1& .129\E1& .126\E2& .694\E2& .487\E2& .562\E2\\
\hline
\text{\cite{CM87}}& .7545& .126& .377\E1& .180\E1& .140\E1& .175\E2& .686\E2& .484\E2& .629\E2\cr  
\hline
\end{array}$   
\end{table}

As is the case in other similar tables, the class number~$1$ appears more
frequently than expected, while most other class groups appear less often.
The data seem to support the heuristic, except that class numbers divisible
by~3 occur considerably less frequently than expected, see the deviations
listed in Table~\ref{tab:dev4}. A similar but even more extreme phenomenon
could be observed in our results in \cite[Tab.~5]{Ma06} for totally real
$\fS_4$-fields, where the predicted value for $h'=9$ was almost double the
observed one.

\begin{table}[htb]
\caption{Deviation of $h'$ from CM-heuristic in percent}   \label{tab:dev4}
$\begin{array}{c|ccccccccccc}
 h'& 1& 3& 5& 7& 9& 11& 13& 15& 17& 19& 21\\
\hline
 \% & 2.5& 8.4& 1.8& .3& 13.6& .2& .1& 15.6& .6& 1.2& 15.8\\
\end{array}$
\end{table}

On the question of whether their heuristic should also apply to the prime~3,
Cohen and Martinet write in the last paragraph of \cite{CM87},
\begin{center}
``about quartic extensions of type $\fA_4$ and $\fS_4$: The prime~3\\
 \emph{could} be bad. However, we think this is not the case."
\end{center}
The data so far seem inconclusive to the author.
Also note the comments of Cohen--Martinet in \cite{CM94}.

\section{Extensions of degree~5 with one real embedding}   \label{sec:S5}
In degree~5 we have enumerated the (necessarily primitive) number fields (up to
isomorphism) with one real embedding and discriminant at most $10^9$.

\subsection{The fields}
The possible Galois groups of the Galois closure are the dihedral group $D_5$
of order~10, the Frobenius group $F_{20}$ of order~20, the alternating group
$\fA_5$, and the symmetric group $\fS_5$.

\begin{prop}   \label{prop:deg5}
 The number fields of degree~$5$ (up to isomorphism) with one real embedding
 and discriminant at most $10^9$ have Galois groups as shown in
 Table~$\ref{tab:gal5}$.
\end{prop}

\begin{table}[htb]
\caption{ Degree~5 fields with $r_1=1$}  \label{tab:gal5}
$\begin{array}{|c|r|r|r|r|}
\hline
    X& N(\fS_5,X)& N(\fA_5,X)& N(F_{20},X)& N(D_5,X) \cr
\hline\hline
       10^4&      69&    .&   .&    2 \cr
       10^5&    1714&    4&   4&   13 \cr
       10^6&   28831&   72&  31&   57 \cr
       10^7&  406505&  460& 123&  196 \cr
       10^8& 5089844& 2335& 411&  664 \cr
       10^9&59505732&10318&1429& 2213 \cr  
\hline
\end{array}$
\end{table}

Some further statistics on the $\fS_5$-fields is given in Table~\ref{tab:S5}.
While the leading term $N(\fS_5,X)\sim B_5\,X$ in the asymptotic has been
proved by Bhargava \cite{Bh10}, to our knowledge there is no prediction yet for
the form of the error term. Defining $\alpha_1,\alpha_2$ as for degree~4 in the
previous section, the
values in Table~\ref{tab:S5} seem to indicate that the error term might be of
the order $O(X^\alpha)$ with $.825\le\alpha\le.876$. In our computation of
totally real $\fS_5$-fields up to $5\cdot 10^9$ we observed that the
corresponding exponent $\alpha$ should lie between $.738$ and $.901$,
consistent with (but weaker than) our bounds here.

Assuming the same form of asymptotic behaviour~(1) as in the case of
$\fS_4$-fields with $B_4$ replaced by $B_5$ and leading error term of the form
$X^{5/6}$, the least square method (with interpolation points at $10^{i/10}$,
$i\ge30$) yields the approximations
$$c'=-1.5121,\qquad c''=.11082,\qquad c'''=1.4331, $$
which predicts 663\,746\,074 $\fS_5$-extensions of discriminant at most $10^{10}$.
The corresponding values of $E(\fS_5,X)$ and $P(\fS_5,X)$ are also printed in
Table~\ref{tab:S5}.

\begin{table}[htb]
\caption{$\fS_5$-fields with $r_1=1$}  \label{tab:S5}
$\begin{array}{|c|r|r|r|r|r|}
\hline
    X& N(\fS_5,X)& E(\fS_5,X)& \qquad\alpha_1& \qquad\alpha_2& P(\fS_5,X)\\ 
\hline\hline
        10^5&     1714&     6921& .768& .924& 1674 \cr
        10^6&    28831&    57519& .793& .913& 288875 \cr
        10^7&   406505&   457000& .809& .898& 405812 \cr
        10^8&  5089844&  3545209& .819& .885& 5091027 \cr
        10^9& 59505732& 26844807& .825& .876& 59507369\cr 
\hline
\end{array}$   
\end{table}

Again, we have counted the discriminant multiplicities for the
$\fS_5$-extensions in our range:

\begin{prop}   \label{prop:mult5}
 The multiplicities at most~$13$ of discriminants of the first $N$
 $\fS_5$-fields with one real embedding, where
 $N\in\{5\cdot 10^6,10^7,5\cdot10^7\}$, are as given in
 Table~$\ref{tab:discs5}$. The discriminant $684\,450\,000=2^4 3^4 5^5 13^2$
 occurs $22$ times, the highest multiplicity in the range.
\end{prop}

\begin{table}[htb] \caption{$k$-tuples of $\fS_5$-discriminants}
  \label{tab:discs5}
$\begin{array}{|r|r|r|r|r|r|r|r|r|r|r|r|r|r|r|}
\hline
 \hbox{\diagbox{$N$}{$k$}}& 1& 2& 3& 4& 5& 6& 7& 8& 9& 10& 11& 12& 13\\
\hline
 5\cdot 10^6& 4245227& 314356& 31604& 5128& 1300& 397& 138& 52& 24& 13& 7& 2& 2\cr
 10^7& 8429602& 649781& 67421& 11125& 2858& 920& 316& 118& 58& 31& 12& 6& 3\cr
 5\cdot 10^7& 41531916& 3460537& 379856& 65550& 16659& 5514& 2075& 885& 375& 179& 92& 34& 23\cr
\hline
\end{array}$  
\end{table}

As already seen in degree~4, the maximal number of non-isomorphic
$\fS_5$-fields with equal discriminant~$D$ grows much slower than $D^{1/2}$.

\subsection{Class groups statistics}
Computation with \gp{} gives:

\begin{prop}   \label{prop:clno5}
 The distribution of $5'$-parts $h'<60$ of the class numbers among the first
 $5\cdot10^7$ $\fS_5$-fields with one real embedding is displayed in
 Table~$\ref{tab:clno5}$; in addition there are~$505$ fields with $h'\ge60$.
\end{prop}

\begin{table}[htb] \caption{$5'$-parts $h'$ of class numbers of $\fS_5$-fields}
  \label{tab:clno5}
$\begin{array}{|c|r||c|r||c|r||c|r||c|r||c|r||c|r||c|r|}
\hline
  h'& n(h')& h'& n(h')& h'& n(h')& h'& n(h')& h'& n& h'& n& h'& n& h'& n\\
\hline
 1& 40617028&  8& 112539& 16& 12066& 23& 2091& 31&  751& 38& 357& 46& 149& 53& 70\\  
 2&  6097592&  9&  58463& 17&  6255& 24& 3345& 32& 1015& 39& 317& 47& 106& 54& 88\\  
 3&  1709540& 11&  26384& 18&  7245& 26& 1652& 33&  687& 41& 201& 48& 218& 56& 78\\  
 4&   902985& 12&  32915& 19&  4344& 27& 1426& 34&  577& 42& 302& 49& 120& 57& 61\\  
 6&   243111& 13&  15005& 21&  3540& 28& 1690& 36&  693& 43& 166& 51&  92& 58& 40\\  
 7&   114393& 14&  15344& 22&  2909& 29&  850& 37&  318& 44& 236& 52& 100& 59& 41\\ 
\hline\end{array}$
\end{table}

In Table~\ref{tab:clno5rel} we compare the relative proportions of $5'$-parts of
class groups to the heuristic predictions of Cohen--Martinet. According to
\cite{CM87} the probability that the $5'$-part $H_K'$ of the class group of an
$\fS_5$-extension $K/\QQ$ with one real embedding is isomorphic to $H$ should
be given by
$$\pr(H_K'=H)=c_5\,\big(|H|^2|\Aut(H)|\big)^{-1}$$
for some explicit constant $c_5=0.7240198...$.

In \cite[Conj.~2.1]{Ma10} we conjectured that a 2-group $H$ of 2-rank $r$
occurs as Sylow 2-subgroup of a class group of an $\fS_5$-field with one real
embedding with probability
\begin{equation}\label{eq:2u=2}
  0.786\ldots\cdot
  \frac{2^{(r^2-r)/2}(2)_{r+2}}{|H|^2\cdot|\Aut(H)|}.
\end{equation}

\begin{table}[ht] \caption{Relative proportions of $5'$-parts $H'$ of $\fS_5$-class groups}
  \label{tab:clno5rel}
$\small\begin{array}{|r|r|r|r|r|r|r|r|r|r|r|r|r|r|r|}
\hline
 \hbox{\diagbox{$N$}{$H'$}}& 1& 2& 3& 4& 2^2& 6& 7& 8& 4\ti2& 2^3& 9& 3^2\cr
\hline
0$--$2\Ex7& .823& .116& .327\E1& .128\E1& .35\E2& .43\E2& .22\E2& .14\E2& .51\E3& .13\E4& .10\E2& .51\E4\\
2$--$4\Ex7& .806& .125& .350\E1& .143\E1& .47\E2& .52\E2& .23\E2& .16\E2& .73\E3& .28\E4& .12\E2& .79\E4\\
4$--$6\Ex7& .802& .128& .357\E1& .148\E1& .51\E2& .55\E2& .24\E2& .17\E2&
.81\E3& .36\E4& .12\E2&  .85\E4\\
\hline
\text{\cite{CM87}}& .724& .181& .402\E1& .227\E1& .75\E2& .10\E1& .25\E2& .28\E2& .14\E2& .67\E4& .15\E2& .19\E3\cr
\text{\cite{Ma10}}& .739& .162& .411\E1& .197\E1& .13\E1& .90\E2& .25\E2& .25\E2& .23\E2& .43\E3& .15\E2& .19\E3\cr
\hline
\end{array}$
\end{table}

As can be seen, the observed values are still quite far from the predictions, a
phenomenon which also occurred for totally real $\fS_5$-fields in
\cite[Tab.~9]{Ma06}. Of course the discrepancy might be due to the very
limited data.

\section{Totally complex $\fS_6$-extensions}   \label{sec:S6}
In degree~6 we have enumerated the primitive totally complex number fields (up
to isomorphism) with absolute discriminant at most $10^9$.

\subsection{The fields}
The only possible primitive Galois group of a totally complex degree~6 field
apart from $\fS_6$ is $\PGL_2(5)\cong\fS_5$.

\begin{prop}   \label{prop:deg6}
 There are $8\,681\,205$ totally complex number fields of degree~$6$ (up to
 isomorphism) with absolute discriminant at most $10^9$ and group $\fS_6$, and
 $1\,108$ further such fields with group $\PGL_2(5)$.
\end{prop}

Some statistics on the $\fS_6$-fields is given in Table~\ref{tab:S6}, with
notation as in the previous sections.

\begin{table}[htb]
\caption{Totally complex $\fS_6$-fields}   \label{tab:S6}
$\begin{array}{|c|r|r|r|r|}
\hline
    X& N(\fS_6,X)& E(\fS_6,X)& \qquad\alpha_1& \qquad\alpha_2\\
\hline\hline
      10^5&      74&    1628& .642& .982  \cr
      10^6&    2276&   14749& .695& .940  \cr
      10^7&   43254&  126999& .729& .927  \cr
      10^8&  654757& 1047773& .753& .910 \cr
      10^9& 8681205& 8344095& .769& .897 \cr
\hline
\end{array}$   
\end{table}

The leading term in the asymptotic is conjectured by Bhargava to be given by
$N(\fS_6,X)\sim B_6\,X$; no prediction for the form of the error term is known
to us. The values in Table~\ref{tab:S6} point towards an error term of the
order $O(X^\alpha)$ with $.769\le\alpha\le.897$. Assuming the same form of
asymptotic expansion~(1) as in the case of $\fS_4$-fields
with $B_4$ replaced by $B_6$ the least square method (with interpolation points
at $10^{i/10}$, $i\ge50$) yields the approximations
$$c'=-.6038,\qquad c''=.068,\qquad c'''=.502,$$
predicting 105\,574\,058 $\fS_6$-extensions of discriminant at most $10^{10}$.

In Table~\ref{tab:discs6} we show statistics on the multiplicity of
discriminants in the range of our computations; no discriminant occurs for
more than~7 distinct fields.

\begin{table}[htb]
\caption{$k$-tuples of $\fS_6$-discriminants}   \label{tab:discs6}
$\begin{array}{|r|r|r|r|r|r|r|r|r|r|}
\hline
 \hbox{\diagbox{$N$}{$k$}} & 1& 2& 3& 4& 5& 6& 7\cr
\hline
        10^6&  968723&  14885&  453&  31&  5&  .& .\cr   
 4\cdot 10^6& 3853637&  68571& 2720& 223& 23&  8& 1\cr   
 8\cdot 10^6& 7686645& 146005& 6242& 544& 68& 15& 2\cr   
\hline
\end{array}$
\end{table}

\subsection{Class groups statistics}

\begin{prop}   \label{prop:clno6}
 The relative proportions of $6'$-parts $H'$ of class groups with $|H'|\le23$
 among the first totally complex $\fS_6$-fields are displayed in
 Table~$\ref{tab:clno6rel}$.
\end{prop}

\begin{table}[ht]
\caption{Relative proportions of $6'$-parts of $\fS_6$-class groups}
  \label{tab:clno6rel}
$\begin{array}{|r|r|r|r|r|r|r|r|r|r|}
\hline
    H'& 1& 5& 7& 11& 13& 17& 19& 23\cr
\hline
0-3\cd10^6& .992& .57\E2& .17\E2& .30\E3& .15\E3& .37\E4& .23\E4& .37\E5\\
3-6\cd10^6& .990& .69\E2& .22\E2& .46\E3& .24\E3& .79\E4& .46\E4& .16\E4\\
6-9\cd10^6& .989& .72\E2& .23\E2& .51\E3& .27\E3& .91\E4& .72\E4&  .24\E4\\
\hline
\text{\cite{CM87}}& .985& .99\E2& .34\E2& .81\E3& .49\E3& .21\E3& .15\E3& .85\E4\\
\hline
\end{array}$
\end{table}

The last line of Table~\ref{tab:clno6rel} shows the predicted probability
$$\pr(H_K'=H)=0.984725...\cdot\big(|H|^2|\Aut(H)|\big)^{-1}$$
from \cite{CM87} that the $6'$-part $H_K'$ of the class group of a totally
complex $\fS_6$-extension equals $H$

\section{Extensions of degree~7 with one real embedding}   \label{sec:S7}
In degree~7 we have enumerated the number fields with one real embedding and
absolute discriminant at most $4\cdot10^9$. Possible Galois groups here are
the dihedral group $D_7$ of order~14, the Frobenius group $F_{42}$ of order~42,
and the full symmetric group.

\begin{prop}   \label{prop:deg7}
 The number fields of degree~$7$ (up to isomorphism) with one real embedding
 and absolute discriminant at most $4\cdot10^9$ have Galois groups as given in
 Table~$\ref{tab:gal7}$.
\end{prop}

\begin{table}[htb]
\caption{ Primitive degree~7 fields with $r_1=1$}  \label{tab:gal7}
$\begin{array}{|c|r|r|r|}
\hline
    X& N(\fS_7,X)& N(F_{42},X)& N(D_7,X) \cr
\hline\hline
        10^6&      392&  .& 1\cr
        10^7&    12698&  .& 3\cr
        10^8&   279385&  2& 12\cr
        10^9&  4690185& 23& 31\cr
 4\cdot 10^9& \ge23611863& 52& 52\cr 
\hline
\end{array}$
\end{table}

Some statistics on the $\fS_7$-fields is given in Table~\ref{tab:S7}. Here,
$E(\fS_7,X)$ and $\alpha_1,\alpha_2$ are as in the previous sections.

\begin{table}[htb]
\caption{$\fS_7$-fields with $r_1=1$}   \label{tab:S7}
$\begin{array}{|c|r|r|r|r|}
\hline
    X& N(\fS_7,X)& E(\fS_7,X)& \qquad\alpha_1& \qquad\alpha_2\\
\hline\hline
        10^6&      392&    17829& .709& .982\cr
        10^7&    12698&   169520& .747& .973\cr
        10^8&   279385&  1542801& .774& .952\cr
        10^9&  4690185& 13531675& .792& .937\cr 
 4\cdot 10^9& \ge23611863& \le49295553& .801& .931\cr 
\hline
\end{array}$ 
\end{table}

The data indicate an error term of the form $O(X^\alpha)$ with
$.801\le\alpha\le .932$. Note, however, that in our range the deviation
$E(\fS_7,X)$ from Bhargava's prediction $B_7\,X$ for $N(\fS_7,X)$ is more than
twice the actual value.

In Table~\ref{tab:discs7} we record the number of $k$-tuples of equal
discriminants among $\fS_7$-fields in our range. There are no
more than nine fields with the same discriminant.

\begin{table}[htb]
\caption{$k$-tuples of $\fS_7$-discriminants}   \label{tab:discs7}
$\begin{array}{|r|r|r|r|r|r|r|r|r|r|r|r|}
\hline
\hbox{\diagbox{$N(\fS_7)$}{$k$}}& 1& 2& 3& 4& 5& 6& 7& 8& 9\cr
\hline
        10^6&  980993&   9129&  224&  14&  3&  1& .& .& .\cr
 5\cdot 10^6& 4876447&  58319& 1975& 192& 28& 10& 2& 1& .\cr 
        10^7& 9724598& 129106& 4793& 504& 99& 39& 7& 1& 1\cr 
\hline
\end{array}$  
\end{table}


\section{Totally complex $\fS_8$-extensions}   \label{sec:S8}
In degree~8 we have enumerated primitive totally complex number fields with
discriminant at most $2\cdot10^{10}$. The possible Galois groups here, apart
from the alternating group $\fA_8$ and the full symmetric group $\fS_8$, are
the transitive groups commonly denoted 8T36, 8T37, 8T43, 8T48. Here 8T36 is the
semidirect product of the elementary abelian group $C_2^3$ with the Frobenius
group $F_{21}$, 8T37 is
$\PSL_2(7)$, 8T43 is $\PGL_2(7)$ and 8T48 is the semidirect product of $C_2^3$
with $\PSL_2(7)$. In our range there is no field with group 8T36 or 8T37. 
Our computations prove the following:

\begin{prop}   \label{prop:min8}
 The number of primitive totally complex number fields of degree~$8$ (up to
 isomorphism) of discriminant at most $2\cdot10^{10}$ is at least as given in
 Table~$\ref{tab:gal8}$.   \par
 The minimal discriminant of octic totally complex fields with group
 $8T43=\PGL_2(7)$ is $418\,195\,493$, and it is $16\,410\,601$ for the
 group $\fA_8$. In either case, there is a unique such field (up to
 isomorphism).
\end{prop}

\begin{table}[htb]
\caption{Primitive totally complex degree~8 fields}  \label{tab:gal8}
$\begin{array}{|c|r|r|r|r|}
\hline
    X& N(\fS_8,X)& N(\fA_8,X)& N(8T48,X)& N(\PGL_2(7),X) \cr
\hline\hline
           10^8&       14261&    9&    4&  . \cr
           10^9&      322794&   92&  102&  1 \cr 
        10^{10}&  \ge5627516&  737&  709& 13 \cr 
 2\cdot 10^{10}& \ge12785743& 1358& 1169& 22\cr 
\hline
\end{array}$
\end{table}

Further data on the $\fS_8$-fields, with notation as in previous sections, are
collected in Table~\ref{tab:S8}. This time, $\alpha_1,\alpha_2$ give only quite
weak information on a possible error term.

\begin{table}[htb]
\caption{Totally complex $\fS_8$-fields}   \label{tab:S8}
$\begin{array}{|c|r|r|r|r|}
\hline
    X& N(\fS_8,X)& E(\fS_8,X)& \qquad\alpha_1& \qquad\alpha_2\\
\hline\hline
           10^7&         373&    24315& .627& .987\cr
           10^8&       14261&   232619& .671& .976\cr
           10^9&      322797&  2146003& .704& .958\cr
        10^{10}&  \ge5627696& \le19060304& .728& .943\cr 
 2\cdot 10^{10}& \ge12792145& \le36583855& .734& .941\cr 
\hline
\end{array}$  
\end{table}

The occurring discriminant multiplicities are shown in Table~\ref{tab:discs8}.
So here at most four fields in the range share the same discriminant.

\begin{table}[htb]
\caption{$k$-tuples of $\fS_8$-discriminants}   \label{tab:discs8}
$\begin{array}{|r|r|r|r|r|r|}
\hline
\hbox{\diagbox{$N(\fS_8)$}{$k$}}& 1& 2& 3& 4\cr
\hline
        10^6&  997556&  1216&   4& .\cr
 5\cdot 10^6& 4982402&  8712&  58& .\cr
        10^7& 9959543& 20009& 145& 1\cr
\hline
\end{array}$  
\end{table}

\section{Extensions of degree~9 with one real embedding}  \label{sec:S9}
In degree~9 we have enumerated primitive number fields with one real
embedding and discriminant at most $4\cdot10^{10}$. Among the possible
Galois groups, only the solvable group 9T16, the group 9T32 isomorphic to
$\PSL_2(8).3$, the alternating group $\fA_9$ and the symmetric group $\fS_9$ do
actually occur (see also \cite{Jo13} for the minimal discriminants for solvable
nonic fields).

\begin{prop}   \label{prop:min9}
 The number of primitive number fields of degree~$9$ (up to isomorphism) with
 one real embedding and discriminant at most $4\cdot10^{10}$ is at least as
 given in Table~$\ref{tab:gal9}$.   \par
 The minimal discriminant of nonic number fields with one real embedding and
 Galois group $\fA_9$ is $92\,371\,321$ and there is a unique field with that
 group and discriminant (up to isomorphism).
\end{prop}

\begin{table}[htb]
\caption{Primitive degree~9 fields with $r_1=1$}   \label{tab:gal9}
$\begin{array}{|c|r|r|r|r|}
\hline
    X& N(\fS_9,X)& N(\fA_9,X)& N(\PSL_2(8).3,X)& N(9T16,X) \cr
\hline\hline
             10^8&      803&   1& .& .\cr
             10^9&    61353&  21& .& .\cr
          10^{10}&  1676194& 196& .& 3\cr
 4\cdot10^{10}& \ge10403964& 764& 3& 9\cr 
\hline
\end{array}$
\end{table}

The results for $\fS_9$-fields are collected in Table~\ref{tab:S9}, with the
usual notation, and the discriminant multiplicities are recorded in
Table~\ref{tab:discs9}. Again, the parameters $\alpha_1,\alpha_2$ are still
quite far apart.

\begin{table}[htb]
\caption{$\fS_9$-fields with $r_1=1$}   \label{tab:S9}
$\begin{array}{|c|r|r|r|r|}
\hline
    X& N(\fS_9,X)& E(\fS_9,X)& \qquad\alpha_1& \qquad\alpha_2\\
\hline\hline
             10^8&      803&   256565& .676& .997\cr
             10^9&    61353&  2512327& .711& .988\cr
          10^{10}&  1676194& 24060606& .738& .977\cr 
 4\cdot10^{10}& \ge10403964& \le92543236& .751& .970\cr 
\hline
\end{array}$  
\end{table}

\begin{table}[htb]
\caption{$k$-tuples of $\fS_9$-discriminants}   \label{tab:discs9}
$\begin{array}{|r|r|r|r|r|r|}
\hline
\hbox{\diagbox{$N(\fS_9)$}{$k$}}& 1& 2& 3& 4\cr
\hline
        10^6&  999094&  453& .& .\cr
 5\cdot 10^6& 4993243& 3374& 3& .\cr 
        10^7& 9983721& 8115& 15& 1\cr 
\hline
\end{array}$  
\end{table}

\section{Totally complex $\fS_{10}$-extensions}   \label{sec:S10}

In degree~10 we have enumerated primitive totally complex number fields of
absolute discriminant at most $4\cdot10^{11}$. The possible Galois groups
are $\PGL_2(9)$, $\Aut(\fS_6)$ and the symmetric group $\fS_{10}$. We did not
find fields with the first two Galois groups in that range.

\begin{prop}   \label{prop:deg10}
 The number of totally complex $\fS_{10}$-fields of degree~$10$ (up to
 isomorphism) with absolute discriminant at most $4\cdot10^{11}$ is at least
 as given in Table~{\rm\ref{tab:S10}}.   \par
 The minimal absolute discriminant for such fields equals $215\,067\,767$, and
 there is a unique such field with that discriminant (up to isomorphism).
\end{prop}

\begin{table}[htb]
\caption{Totally complex $\fS_{10}$-fields}   \label{tab:S10}
$\begin{array}{|c|r|r|r|r|}
\hline
    X& N(\fS_{10},X)& E(\fS_{10},X)& \qquad\alpha_1& \qquad\alpha_2\\
\hline\hline
          10^9&         817&      267583& .603& .997\cr
       10^{10}&       58671&     2625329& .642& .989\cr 
       10^{11}&  \ge1634855& \le25205145& .673& .979\cr 
 4\cdot10^{11}& \ge10030338& \le97329662& .689& .974\cr 
\hline
\end{array}$  
\end{table}

The count of multiple discriminants is shown in Table~\ref{tab:discs10}.

\begin{table}[htb]
\caption{$k$-tuples of $\fS_{10}$-discriminants}   \label{tab:discs10}
$\begin{array}{|r|r|r|r|r|r|}
\hline
\hbox{\diagbox{$N(\fS_{10})$}{$k$}}& 1& 2& 3\cr
\hline
        10^6&  999902&  49& .\cr
 5\cdot 10^6& 4999300& 350& .\cr
        10^7& 9998370& 815& .\cr
\hline
\end{array}$  
\end{table}

\section{Extensions of degree~11 with one real embedding}    \label{sec:S11}
Finally, in degree~11 we have enumerated number fields with one real embedding
and discriminant at most $2\cdot10^{12}$. Possible Galois groups are the
dihedral group $D_{11}$, the Frobenius group $F_{110}$ and the symmetric group
$\fS_{11}$. There is one $D_{11}$-extension and no $F_{110}$-extension in that
range. The data for the $\fS_{11}$-extensions are as follows:

\begin{prop}   \label{prop:min11}
 The number of $\fS_{11}$-fields of degree~$11$ (up to isomorphism) with one
 real embedding and absolute discriminant at most $2\cdot10^{12}$ is at least
 as given in Table~$\ref{tab:S11}$.\par
 The minimal absolute discriminant of such fields equals $5\,781\,612\,911$,
 and there is a unique field with that discriminant (up to isomorphism).
\end{prop}

\begin{table}[htb]
\caption{$\fS_{11}$-fields with $r_1=1$}   \label{tab:S11}
$\begin{array}{|c|r|r|r|r|}
\hline
    X& N(\fS_{11},X)& E(\fS_{11},X)& \qquad\alpha_1& \qquad\alpha_2\\
\hline\hline
        10^{10}&         789&      2747011& .644& .9996\cr
        10^{11}&      162093&     27315907& .676& .996\cr 
        10^{12}&  \ge6698503& \le268081497& .702& .990\cr 
 2\cdot 10^{12}& \ge16863137& \le532696863& .709& .991\cr 
\hline
\end{array}$  
\end{table}

The value of $\alpha_2$ in the last line of Table~\ref{tab:S11} hints that our
enumeration in this case seems far from complete, yet.

\begin{table}[htb]
\caption{$k$-tuples of $\fS_{11}$-discriminants}   \label{tab:discs11}
$\begin{array}{|r|r|r|r|r|r|}
\hline
   \quad k& 1& 2& 3\cr
   N(\fS_{11})\quad& & & \cr
\hline
        10^6&  999978& 11& .\cr
 5\cdot 10^6& 4999830& 85& .\cr
        10^7& 9999600& 200& .\cr
\hline
\end{array}$  
\end{table}

\section{Totally real $\fS_n$-number fields}   \label{sec:tr}

We have also continued our computation from \cite{Ma06} of totally real
primitive number fields and extended them to degrees up to~11. Here we also
used the further restrictions in Hunter's algorithm as described in \cite{Ma06}.

The lists are as yet incomplete in degree $n\ge6$; in particular in
degrees $n\ge9$ we expect many more fields than encountered so far. The number
of fields obtained are listed in Table~\ref{tab:nr-tr}, the qualitative data
are displayed in Figure~\ref{tab:plot-tr}, with the same understanding as for
Figure~\ref{tab:plot}, that is, the plotted lines indicate the number of
(computed) totally real $\fS_n$-extensions of $\QQ$, with $n$ increasing from~4
at the left to~11 at the right.

\begin{table}[htb]
\caption{Number of totally real $\fS_n$-extensions}   \label{tab:nr-tr}
$$\begin{array}{c|rrc}
 n& X& N_{tr}(\fS_n,X)\\
\hline
  4&           10^9&     17,895,702\\  
  5&    5\cdot 10^9&     14,153,774\\  
  6& 5\cdot 10^{10}& \ge 17,323,021\\  
  7& 4\cdot 10^{11}& \ge 11,180,749\\  
  8& 5\cdot 10^{12}&  \ge 9,108,220\\  
  9&        10^{14}&  \ge 8,549,839\\  
 10& 3\cdot 10^{15}&  \ge 8,820,242\\  
 11&        10^{17}&  \ge 5,119,912\\  
\end{array}$$
\end{table}

\begin{table}[htbp]
\caption{Totally real $\fS_n$-number fields, $4\le n\le 11$}   \label{tab:plot-tr}
\begin{tikzpicture}
\begin{loglogaxis}[
    legend style = { at = {(0.25,0.95)}},
    xmin = 500000, xmax = 2*10^17,    ymin = 400, ymax = 2*10^7,
    xtick distance = 10,    ytick distance = 10,
    grid = none,    minor tick num = 1,
    width = \textwidth,    height = 15cm, 
    xlabel = {$X$},    ylabel = {$N(\fS_n,X)$},]
 
\addplot[domain=3*10^7:7*10^8, dashed]{0.025347714*x}; 
\addplot[ thin, smooth ]
coordinates {
(39810.7170,118)
(63095.7344,225)
(100000.000,449)
(158489.319,852)
(251188.643,1524)
(398107.170,2697)
(630957.344,4793)
(1000000.00,8301)
(1584893.19,14398)
(2511886.43,24691)
(3981071.70,42244)
(6309573.44,71807)
(10000000.0,120622)
(15848931.9,202514)
(25118864.3,337547)
(39810717.0,560622)
(63095734.4,927683)
(100000000,1529634)
(158489319,2515160)
(251188643,4122288)
(398107170,6739945)
(630957344,10994651)
(1000000000,17895702) };  
\addplot[domain=10^8:2.7*10^9, dashed]{0.005756702*x}; 
\addplot[ thin, smooth ]
coordinates {
(251188.643,52)
(398107.170,104)
(630957.344,213)
(1000000.00,409)
(1584893.19,765)
(2511886.43,1499)
(3981071.70,2820)
(6309573.44,5165)
(10000000.0,9461)
(15848931.9,16982)
(25118864.3,30135)
(39810717.0,53260)
(63095734.4,92887)
(100000000.,162022)
(158489319.,279749)
(251188643.,480115)
(398107170.,818014)
(630957344.,1387546)
(1000000000,2341055)
(1584893192,3935609)
(2511886431,6593530)
(3981071705,11000768)
(5000000000,14153774) };  
\addplot[domain=5*10^8:1.3*10^10, dashed]{0.00113502007*x}; 
\addplot[ thin, smooth ]
coordinates {
(6309573.44,77)
(10000000.0,177)
(15848931.9,370)
(25118864.3,782)
(39810717.0,1608)
(63095734.4,3288)
(100000000.,6513)
(158489319,12512)
(251188643,23524)
(398107170,43858)
(630957344,80961)
(1000000000,147553)
(1584893192,265702)
(2511886431,475090)
(3981071705,842132)
(6309573444,1479961)
(10000000000,2585968)
(15848931924,4489912)
(25118864315,7744036)
(39810717055,13283501)
(50000000000,17317471) };  
\addplot[domain=2*10^9:9*10^10, dashed]{0.00017354146*x }; 
\addplot[ thin, smooth ]
coordinates {
(100000000.0,71)
(158489319.2,173)
(251188643.1,379)
(398107170.5,803)
(630957344.4,1673)
(1000000000.0,3427)
(1584893192.4,6871)
(2511886431.5,13539)
(3981071705.5,26368)
(6309573444.8,50839)
(10000000000,96299)
(15848931924,180146)
(25118864315,333174)
(39810717055,610703)
(63095734448,1108343)
(100000000000,1996894)
(158489319246,3569139)
(251188643151,6332229)
(400000000000,11180749) };  
\addplot[domain=10^10:6*10^11, dashed]{0.0000235124116*x }; 
\addplot[ thin, smooth ]
coordinates {
(3981071705.535,164)
(6309573444.802,402)
(10000000000.00,927)
(15848931924.61,2157)
(25118864315.10,4659)
(39810717055.35,10014)
(63095734448.02,20666)
(100000000000.0,42065)
(158489319246.1,83573)
(251188643151.0,163277)
(398107170553.5,315597)
(630957344480.2,600427)
(1000000000000,1132306)
(1584893192461,2108530)
(2511886431510,3863925)
(3981071705535,6910616)
(5000000000000,9108220) };  
\addplot[domain=10^11:5*10^12, dashed]{0.00000272347149294*x }; 
\addplot[ thin, smooth ]
coordinates {
(63095734448.0,57)
(100000000000.0,173)
(158489319246.1,421)
(251188643150.9,1131)
(398107170553.4,2690)
(630957344480.1,6095)
(1000000000000,13329)
(1584893192461,28363)
(2511886431510,58982)
(3981071705535,120287)
(6309573444802,240991)
(10000000000000,474970)
(15848931924610,911437)
(25118864315100,1694700)
(39810717055350,3031373)
(63095734448020,5198401)
(100000000000000,8549839) };  
\addplot[domain=2*10^12:5*10^13, dashed]{0.00000028401716*x }; 
\addplot[ thin, smooth ]
coordinates {
(2511886431509.58,147)
(3981071705534.97,386)
(6309573444801.93,960)
(10000000000000.0,2388)
(15848931924611.1,5748)
(25118864315095.8,13218)
(39810717055349.7,29517)
(63095734448019.3,63896)
(100000000000000.,134696)
(158489319246111.,272235)
(251188643150958.,524141)
(398107170553497.,958965)
(630957344480193.,1670337)
(1000000000000000,2780936)
(1584893192461000,4444343)
(2511886431510000,6861260)
(3000000000000000,8049303) };  
\addplot[domain=4*10^13:5*10^14, dashed]{.0000000264335054*x }; 
\addplot[ thin, smooth ]
coordinates {
(100000000000000.0,133)
(158489319246111.3,474)
(251188643150958.0,1507)
(398107170553497,4089)
(630957344480193,10375)
(1000000000000000,24235)
(1584893192461113,52593)
(2511886431509580,106752)
(3981071705534972,202333)
(6309573444801932,362207)
(10000000000000000,616607)
(15848931924611135,1001710)
(25118864315095801,1571456)
(39810717055349725,2390726)
(63095734448019325,3539162)
(100000000000000000,5119912) };  
\end{loglogaxis}
 \end{tikzpicture}
\end{table}

\newpage

\end{document}